\documentclass[12pt,a4paper]{article}
\usepackage[sumlimits]{amsmath}
\usepackage{amsthm}
\usepackage{amssymb}
\usepackage{amsfonts}
\usepackage{graphicx}
\usepackage[T1]{fontenc}
\usepackage{enumerate}

\newcommand{\fsl}{\mathfrak{sl}}
\newcommand{\frakg}{\mathfrak{g}}

\newcommand{\bbR}{\mathbb{R}}

\renewcommand{\phi}{\varphi}

\newcommand{\set}[1]{\ensuremath{\left\{#1\right\}}}

\newtheorem{proposition}{Proposition}
\newtheorem{lemma}[proposition]{Lemma}
\newtheorem{theorem}[proposition]{Theorem}
\theoremstyle{definition}

\def\clap#1{\hbox to 0pt{\hss#1\hss}}

\def\mathrlap{\mathpalette\mathrlapinternal}

\def\mathrlapinternal#1#2{%
\rlap{$\mathsurround=0pt#1{#2}$}}

\begin{document}

\title{Bilinear forms on $\mathfrak{sl}_2$-modules and a hypergeometric identity}

\author{Johan K\aa hrstr\"om}%

\maketitle

\abstract{In this paper we study properties of a
	certain bilinear form on finite dimensional $\fsl_2(\bbR)$-modules,
	and how these properties behave with respect to tensor products of
	modules. An attempt to determine the signature of this form
	leads to an interesting identity for the hypergeometric series
	${}_3F_2$, which is known as the Karlsson-Minton identity.}

\section{Introduction}

Let $X, Y, H$ be the standard basis of the Lie-algebra $\frakg=\fsl_2(\bbR)$.
These elements satisfy the following defining relations:
\begin{equation}
	\label{eq:lierel}
	[X, Y] = H,\qquad [H, X] = 2X,\qquad [H, Y] = -2Y.
\end{equation}
Furthermore, let $*$ be the anti-involution on $\frakg$ given by
\begin{equation}
	\label{eq:antiinv}
	X^* = X,\qquad Y^* = Y,\qquad H^* = -H.
\end{equation}

Let $V$ be a finite dimensional vector space over $\bbR$,
and $Q$ be a non-degenerate symmetric bilinear form on $V$.
A \emph{$*$-representation} of $\frakg$ is a
representation $\phi:\fsl_2(\bbR)\rightarrow \mathfrak{gl}(V)$ such that
\[
	Q\bigl(\phi(x)u, v\bigr) = Q\bigl(u, \phi(x^*)v\bigr)
\]
for all $x\in\frakg$, $u, v\in V$.
We will say that $V$ is a $*$-representation of $\frakg$ with respect
to $Q$, and write $xv$ instead of $\phi(x)v$.
Given two $*$-re\-presentations $U$ and $V$, with respect to the bilinear
forms $Q$ and $R$ respectively, the tensor product $U\otimes V$ is also a $*$-representation,
with respect to the induced bilinear form $Q\otimes R$.

For $m\in\{0, 1, 2, \dots\}$, denote by $V_m$ the $(m+1)$-dimensional irreducible
representation of $\frakg$, which has the basis $\{e_{-m+2i}\mid 0\leq i\leq m\}$,
and where the action of our generators in this basis is given by the following formulae
(for simplicity, we assume $e_{-m-2} = e_{m+2} =0$), see \cite[Section 11]{fh},
\begin{equation}
	\label{eq:gaction}
	\begin{split}
	He_{m-2i} &= ie_{m-2i},\\
	Ye_{m-2i} &= e_{m-2(i+1)},\\
	Xe_{m-2i} &= i(m-i+1)e_{m-2(i-1)}.
	\end{split}
\end{equation}

Let $m, n\in\{0, 1, 2, \dots\}$. Let further $Q$ and $R$ be non-degenerate
symmetric bilinear forms on $V_m$ and $V_n$ respectively, which make
these representations of $\frakg$ into $*$-representations.

For $0\leq k\leq\min\{m, n\}$, let $P_{-m-n+2k}$ denote the kernel of $Y$ when acting
on the weight subspace of $V_m\otimes V_n$ of weight $-m-n+2k$ (which is
by definition the set of all $v\in V_m\otimes V_n$ such that
$Hv = (-m-n+2k)v$).
The form $Q\otimes R$ induces the following symmetric bilinear form on
$P_{-m-n+2k}$:
\begin{equation}
	\label{eq:omega}
	\begin{array}{ccc}
	\omega_k : P_{-m-n+2k}\times P_{-m-n+2k} & \rightarrow &\bbR\\
	(v, w)& \mapsto &Q\otimes R(v, X^{m+n-2k}w)
	\end{array}
\end{equation}
The following question was formulated by K.-H.~Fieseler in connection to
the study of Hodge-Riemann relations for polytopes, see \cite{khf}.

\medskip 

\noindent{\bf Question:} {\it Are the forms $\omega_k$
positive/negative definite, alternating with $k$, $0\leq k\leq\min\{m,n\}$?}

\medskip

In \cite{khf}, this question is answered using some abstract arguments. In this
paper, we propose a brute-force solution based on explicit calculation.
The advantage of the latter method is that it leads to an
identity for the hypergeometric series ${}_3F_2$, which is known as
the Karlsson-Minton identity \cite{mint, karl}.
We refer the reader to \cite{pwz} for a survey on
hypergeometric identities. It is interesting to see such an identity
arise from representation theory of Lie algebras.

This is the second draft of this paper, and when writing the
first draft I was not aware of the previous occurence of the
identity in \cite{mint, karl}. I would like to thank
Hjalmar Rosengren and Michael Schlosser
for pointing this out to me.

In Section~2, we show how the positive answer to the above question reduces
to a ${}_3F_2$ identity. This essentially amounts to determining the
sign of certain Clebsch-Gordan coefficients of the tensor product
$V_m\otimes V_n$. For a more general treatment of determining
Clebsch-Gordan coefficiets, see~\cite[Chapter 8]{vk}.

\section{Solution of the problem}

Let $m\in\set{0, 1, 2,\dots}$ and $Q:V_m\times V_m\rightarrow\bbR$ be
a non-degenerate symmetric bilinear form such that $V_m$ becomes a $*$-representation
with respect to $Q$.

\begin{lemma}
	\label{lem:qform}
	The bilinear form $Q$ satisfies the following:
	\begin{enumerate}[(i)]
		\item \label{lem:qform.1}
			$Q(e_i, e_j)=0$ for all $i\neq -j\in\set{-m, -m+2, \dots, m-2, m}$.
		\item \label{lem:qform.2}
			$Q(e_i, e_{-i}) = Q(e_j, e_{-j})$ for all $i, j\in\set{-m, -m+2, \dots, m-2, m}$.
	\end{enumerate}
\end{lemma}

\begin{proof}
	From \eqref{eq:antiinv} we get that
\[
	iQ(e_i, e_j) = Q(He_i, e_j) = Q(e_i, -He_j) = -jQ(e_i, e_j)
\]
for all $i, j\in\set{-m, -m+2, \dots, m-2, m}$, so $Q(e_i, e_j)=0$ provided that
$i\neq-j$. This proves \eqref{lem:qform.1}.

Furthermore,
\[
	Q(e_i, e_{-i}) = Q(Ye_{i+2}, e_{-i}) = Q(e_{i+2}, Ye_{-i}) = Q(e_{i+2}, e_{-i-2})
\]
for all $i\in\set{-m, -m+2, \dots, m-2}$. This proves \eqref{lem:qform.2}.

\end{proof}

Recall from \cite[Section 11]{fh} that for $m, n\in\set{0, 1, 2, \dots}$ we have
\begin{equation}
	\label{eq:kg}
	V_m\otimes V_n = V_{\lvert m-n\rvert}\oplus V_{\lvert m-n\rvert+2}\oplus\cdots\oplus V_{m+n}.
\end{equation}
In particular, we can depict the explicit basis and subspaces $P_i$'s in the module
$V_m\otimes V_n$ as follows:

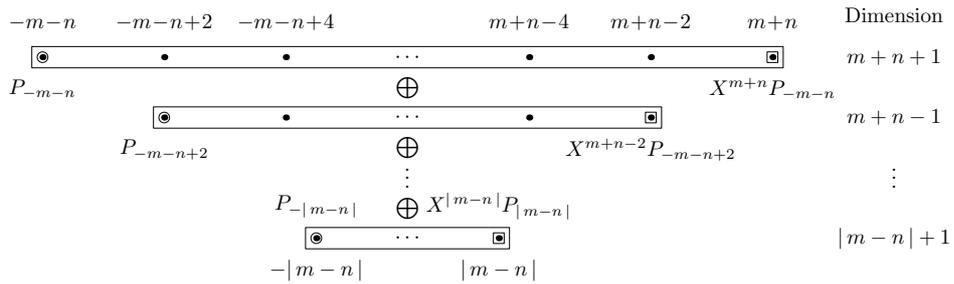
\begin{center}
	\includegraphics[scale=0.8]{irred.1}
\end{center}

To distinguish the bases of $V_m$ and $V_n$ we assume that the set $\{e_{-m}, e_{-m+2}, \dots, e_m\}$
is a basis of $V_m$ and $\set{\tilde e_{-n}, \tilde e_{-n+2}, \dots, \tilde e_n}$ is
a basis of $V_n$. Let $Q$ and $R$ be symmetric bilinear forms on $V_m$ and
$V_n$ respectively, which make these representations into $*$-representations.
By Lemma~\ref{lem:qform}, there exist $q, r\in\bbR$ such that $Q(e_i, e_{-i})=q$
and $R(\tilde e_i, \tilde e_{-i})=r$ for all appropriate $i$.

The subspace of $V_m\otimes V_n$ of weight $-m-n+2k$ is spanned by the vectors
\begin{equation}
	\label{eq:pmnkspan}
	e_{-m}\otimes\tilde e_{-n+2k},\; e_{-m+2}\otimes\tilde e_{-n+2(k-1)},\;\dots,\;
	e_{-m+2k}\otimes\tilde e_{-n}.
\end{equation}

\begin{lemma}
	\label{lem:pbasis}
	For each $0\leq k\leq\min\set{m, n}$ the vector
	\[
	b_{-m-n+2k} = \sum_{i=0}^k(-1)^ie_{-m+2i}\otimes\tilde e_{-n+2(k-i)}.
	\]
	forms a basis of the vector space $P_{-m-n+2k}$.
\end{lemma}

\begin{proof}
	From the decomposition \eqref{eq:kg} it follows that $\dim P_{-m-n+2k}=1$.
	Hence, to find a basis of $P_{-m-n+2k}$ it is sufficient to
	find a non-zero element of weight $-m-n+2k$ in $V_m\otimes V_n$ annihilated
	by $Y$. As the vectors from \eqref{eq:pmnkspan} form a basis of the $(-m-n+2k)$-weight subspace
	of $V_m\otimes V_n$, the latter problem reduces to the solution of the
	following system of linear equations with unknowns $a_i$'s:
	\begin{align*}
 		0 &= Y\Bigl(\sum_{i=0}^ka_i\,e_{-m+2i}\otimes\tilde e_{-n+2(k-i)}\Bigr) \\
		&= \sum_{i=0}^ka_i\,\bigl(Ye_{-m+2i}\otimes\tilde e_{-n+2(k-i)} + e_{-m+2i}\otimes Y\tilde e_{-n+2(k-i)}\bigr) \\
 		&= \sum_{i=1}^{k-1}\Bigl(a_i\,e_{-m+2(i-1)}\otimes\tilde e_{-n+2(k-i)} + a_i\,e_{-m+2i}\otimes\tilde e_{-n+2(k-i-1)}\Bigr) \\
 		&\qquad +a_0\,e_{-m}\otimes\tilde e_{-n+2(k-1)} + a_k\,e_{-m+2(k-1)}\otimes\tilde e_{-n} \\
 		&= \sum_{i=0}^{k-1}(a_{i+1}+a_i)\,e_{-m+2i}\otimes\tilde e_{-n+2(k-i-1)}.
	\end{align*}
	This implies that $a_i=-a_{i+1}$ for all $i$, and the claim follows.
\end{proof}

To proceed we will need the following key result.

\begin{lemma}
	\label{lem:ouridentity}
	For all integers $k, l, m, n$ with $0\leq l\leq k\leq\min\set{m, n}$, the following
	identity holds
	\[
		\sum_{i=\max(0,2k-l-n)}^{\min(k, m-l)}\!\!\!\!\!\!\!\!\!\!\!\!\!
		(-1)^i\frac{(m\!-\!i)!\,(n\!-\!k\!+\!i)!}{i!\,(k\!-\!i)!\,(m\!-\!l\!-\!i)!\,(n\!+\!l\!-\!2k\!+\!i)!}
		= (-1)^{k+l}.
	\]
\end{lemma}

\begin{proof}
	This identity can be restated as the Karlssor-Minton indentity for
	the hypergeometric series ${}_3F_2$. For a proof of the Karlsson-Minton
	identity, see~\cite{mint}.
\end{proof}

To simplify the notation, denote $\mathfrak{s}_k = m+n-2k$.

\begin{lemma}
	\label{lem:xsbs}
	\[
		X^{\mathfrak{s}_k}b_{-\mathfrak{s}_k} =
		\mathfrak{s}_k!\sum_{l=0}^k(-1)^{k+l}\frac{(m-l)!\,(n-k+l)!}{l!\,(k-l)!}e_{m-2l}\otimes\tilde e_{n-2(k-l)}.
	\]
\end{lemma}

\begin{proof}
	Direct calculation yields
\begin{align*}
	X^{\mathfrak{s}_k}\,b_{-\mathfrak{s}_k}
	&=	X^{\mathfrak{s}_k}\Bigl(\sum_{i=0}^k(-1)^ie_{-m+2i}\otimes\tilde e_{-n+2(k-i)}\Bigr) \\
	&= \sum_{j=0}^{n+m-2k}\!\sum_{i=0}^k(-1)^i\binom{m\!+\!n\!-\!2k}{j}
		(X^je_{-m+2i})\otimes\Bigl(X^{\raisebox{1ex}{$\scriptstyle\mathrlap{n+m-2k-j}$}}\,\tilde e_{-n+2(k-i)}\Bigr) \\
	&= \sum_{j=0}^m\,\sum_{i=m-k-j}^{m-j}(-1)^i\frac{(m+n-2k)!}{j!\,(m+n-2k-j)!}\times \\
	&\qquad\times\frac{(m-i)!\,(j+i)!}{(m-i-j)!\,i!}\cdot\frac{(n-k+i)!\,(n+m-k-i-j)!}{(k+i+j-m)!\,(k-i)!}\times\\
	&\qquad\qquad\times e_{-m+2(i+j)}\otimes\tilde e_{-n+2(n+m-k-i-j)} \\
	&= (\;\text{substitution: }l=m-i-j\;) \\
	&= \mathfrak{s}_k!\sum_{l=0}^k\frac{(m-l)!\,(n-k+l)!}{l!\,(k-l)!}\times\\
	&\qquad\times\Biggl(\sum_{i=\max\{0, 2k-n-l\}}^{\min\{k,m-l\}}
		\!\!\!\!\!\!\!\!\!(-1)^i\frac{(m-i)!\,(n-k+i)!}{i!\,(k-i)!\,(m-l-i)!\,(n+l-2k+i)!}\Biggr) \times\\
	&\qquad\qquad\times e_{m-2l}\otimes\tilde e_{n-2(k-l)} \\
	&=	\mathfrak{s}_k!\sum_{l=0}^k(-1)^{k+l}\frac{(m-l)!\,(n-k+l)!}{l!\,(k-l)!}e_{m-2l}\otimes\tilde e_{n-2(k-l)},
\end{align*}
	where for the last equality we use Lemma~\ref{lem:ouridentity}.
\end{proof}

We are now ready to answer the Question from the introduction.

\begin{theorem}
	The forms $\omega_k$ are positive/negative, alternating with
	$k$, $0\leq k\leq\min\set{m, n}$.
\end{theorem}

\begin{proof}
	Since $P_{-\mathfrak{s}_k}$ is spanned by $b_{-\mathfrak{s}_k}$, it suffices to verify that
	the value
	$\omega_k(b_{-\mathfrak{s}_k}, b_{-\mathfrak{s}_k})$
	is positive/negative, alternating with $k$.
\begin{align*}
	\omega_k(b_{-\mathfrak{s}_k}, b_{-\mathfrak{s}_k}) &=
	(Q\otimes R)(b_{-\mathfrak{s}_k}, X^{\mathfrak{s}_k}\,b_{-\mathfrak{s}_k}) \\
	&= \mathfrak{s}_k!\sum_{i=0}^k\sum_{l=0}^k(-1)^{i+k+l}\frac{(m-l)!\,(n-k+l)!}{l!\,(k-l)!}\times\\
	&\qquad\times Q(e_{-m+2i},e_{m-2l})R(\tilde e_{-n+2(k-i)},\tilde e_{n-2(k-l)}) \\
	&= \mathfrak{s}_k!\,(-1)^kqr\sum_{l=0}^k\frac{(m-l)!\,(n-k+l)!}{l!\,(k-l)!},
\end{align*}
	and the claim follows.
\end{proof}

\bigskip

\noindent{\bf Acknowledgements.}
I would like to thank Volodymyr~Mazorchuk for all his help, and
Svante~Janson, Tobias~Ekholm, Herbert~Wilf, Doron~Zeilberger, Christian~Krattenthaler
and Ganna~Kudryavtseva for providing proofs for the identity
of Lemma~\ref{lem:ouridentity}. I would also like to thank
Hjalmar~Rosengren and Michael~Schlosser for pointing out that
this identity is in fact the Karlsson-Minton identity
for ${}_3F_2$. Finally, I would like to thank Anatoliy~Klimyk
for pointing out the connection to the problem
of determining Clebsch-Gordan coefficients.

\vspace{1cm}

\noindent
Department of Mathematics, Uppsala University, Box 480, SE-75106 Uppsala, Sweden,
e-mail, {\tt johan.kahrstrom@math.uu.se}.

\end{document}